\documentclass[12pt]{amsart}
\usepackage{amscd,amssymb,amsthm,verbatim}
\newcommand{\dvol}{\operatorname{dvol}}

\newcommand{\Hess}{\operatorname{Hess}}

\newcommand{\R}{{\mathbb R}}

\newcommand{\Ric}{\operatorname{Ric}}
\newcommand{\Riem}{\operatorname{Riem}}

\newcommand{\Tr}{\operatorname{Tr}}

\newcommand{\vol}{\operatorname{vol}}
\newcommand{\Z}{{\mathbb Z}}

\numberwithin{equation}{section}
\theoremstyle{plain}

\newtheorem{theorem}[equation]{Theorem}

\theoremstyle{remark}

\usepackage{graphicx}
\input{amssym.def}
\input{amssym}
\input{epsf}
\usepackage{a4wide}

\def\R{\mathbb R}

\def\Z{\mathbb Z}

\def\cal{\mathcal}

\def\Hess{\mathop{{\rm Hess}\,}}

\def\2dr#1#2{\left. \frac{d^2}{d{#1}^2} \right |_{#2}}
\def\d2#1{\frac{d^2}{d{#1}^2}}

\def\begeq{\begin{equation}}
\def\endeq{\end{equation}}
\def\begar{\begin{eqnarray}}
\def\endar{\end{eqnarray}}
\def\begar*{\begin{eqnarray*}}
\def\endar*{\end{eqnarray*}}
\def\begal{\begin{align}}
\def\endal{\end{align}}
\def\begal*{\begin{align*}}
\def\endal*{\end{align*}}

\theoremstyle{definition}

\theoremstyle{remark}

\newtheorem*{Thm*}{Theorem}
\newtheorem*{Lem*}{Lemma}
\newtheorem*{Conj*}{Conjecture}
\newtheorem*{Cor*}{Corollary}
\newtheorem*{Def*}{Definition}
\newtheorem*{Prop*}{Proposition}
\newtheorem*{Exo*}{Exercise}
\newtheorem*{Exs*}{Examples}
\newtheorem*{Ex*}{Example}
\newtheorem*{Rk*}{Remark}
\newtheorem*{Rks*}{Remarks}

\begin{document}
\title{Remark about Scalar Curvature and Riemannian Submersions}
\author{John Lott}
\address{Department of Mathematics\\
University of Michigan\\
Ann Arbor, MI  48109-1043\\
USA}
\email{lott@umich.edu}
\thanks{Research supported by NSF grant DMS-0306242 and the
Miller Institute}
\date{May 10, 2005}
\begin{abstract}
We consider modified scalar curvature functions for Riemannian manifolds 
equipped with smooth measures.  Given a Riemannian submersion whose
fiber transport is measure-preserving up to constants, we 
show that the modified scalar curvature of the base is
bounded below in terms of the scalar curvatures of the total space and
fibers. We give an application concerning scalar curvatures of smooth
limit spaces arising in bounded curvature collapses.
\end{abstract}
\maketitle

\section{Introduction} 

An important result in the study of manifolds of nonnegative
sectional curvature is O'Neill's theorem, which says that
sectional curvature is nondecreasing under a Riemannian
submersion
\cite[Chapter 9]{Besse (1987)}. In a previous paper we 
showed that there is a Ricci
analog of O'Neill's theorem, provided that one uses a modified
Ricci tensor
and assumes that the fiber transport of the Riemannian submersion
preserves measures up to multiplicative constants
\cite[Theorem 2]{Lott (2003)}. We used this to show that in the
case of a bounded curvature collapse with a smooth limit space,
a lower Ricci bound on the collapsing manifolds implies a lower
bound on the modified Ricci curvature of the limit space
\cite[Theorem 3]{Lott (2003)}.

These results can be viewed as precursors to a synthetic treatment of
Ricci curvature for metric-measure spaces 
\cite{Lott-Villani (2004),Sturm (2004)}. In particular,
\cite[Theorem 2]{Lott (2003)} is given a synthetic
proof in \cite[Corollary 7.52]{Lott-Villani (2004)}, and
\cite[Theorem 3]{Lott (2003)} follows from
\cite[Corollary 7.45]{Lott-Villani (2004)}.
(To be precise, with the notation of \cite{Lott-Villani (2004)},
if $N < \infty$ then one must restrict to the case of 
nonnegative Ricci curvature.)

In the present note we consider scalar curvature analogs of  
\cite[Theorems 2 and 3]{Lott (2003)}.
A modified scalar curvature function for a smooth Riemannian manifold,
equipped with a smooth measure, was defined by Perelman
\cite[Section 1.3]{Perelman (2002)}.
Let $(M, g)$ be a Riemannian manifold. Let
$\dvol_M$ denote the Riemannian density on $M$.  Let $\phi$ be a smooth
positive function on $M$.  Then $(M, \phi \: \dvol_M)$ is
a smooth metric-measure space. Perelman's modified scalar curvature is
\begin{equation} \label{1.1}
R_\infty \: = \: R \: - \: 2 \: \frac{\nabla^2 \phi}{\phi} \: + \:
\frac{|\nabla \phi|^2}{\phi^2}.
\end{equation}
It is designed so that a weighted Lichnerowicz identity holds. In
particular, if $M$ is a closed spin manifold and $R_\infty \: > \: 0$
then $\widehat{A}(M) \: = \: 0$. 
Perelman's entropy functional is the corresponding total scalar
curvature, i.e. ${\cal F} \: = \: \int_M R_\infty \: \phi \: \dvol_M$.
Note that $R_\infty$ is {\em not} the
trace of the Bakry-\'Emery tensor
$\Ric_\infty \: = \: \Ric \: - \: \frac{\Hess \phi}{\phi} \: + \:
\frac{d \phi}{\phi} \otimes \frac{d \phi}{\phi}$.

More generally, for $q \in (0, \infty)$, we define
\begin{equation} \label{1.2}
R_q \: = \: R \: - \: 2 \: \frac{\nabla^2 \phi}{\phi} \: + \:
\left( 1 \: - \: \frac{1}{q} \right) \: \frac{|\nabla \phi|^2}{\phi^2}.
\end{equation}

We now state a scalar curvature 
analog of the Ricci O'Neill theorem of
\cite[Theorem 2]{Lott (2003)}.
Suppose that a Riemannian submersion  $p \: : \: M \rightarrow B$
has compact fiber $F$. Put $F_b \: = \: p^{-1}(b)$. Given
a smooth curve $\gamma \: : \: [0,1] \rightarrow B$ and a point
$m \in F_{\gamma(0)}$, let $\rho(m)$ be the endpoint 
$\overline{\gamma}(1)$ of the horizontal lift
$\overline{\gamma}$ of $\gamma$ that
starts at $m$. Then $\rho$ is the
fiber transport diffeomorphism
from $F_{\gamma(0)}$ to $F_{\gamma(1)}$. 

Given the positive function $\phi^M$ on $M$, define
$\phi^B$, a smooth positive function on $B$, by
\begin{equation} \label{1.3} 
p_*(\phi^M \: \dvol_M) \: = \: \phi^B \: \dvol_B.
\end{equation}
Let $\phi^F$ denote the restriction of $\phi^M$ to a fiber $F$.
Let $R^M_\infty$, $R^B_\infty$, $R^B_q$ and $R^F_\infty$
denote the modified scalar curvature functions.
Let $\dvol_F$ denote the fiberwise Riemannian density. Let
$\int_F$ denote fiberwise integration.

\begin{theorem} \label{theorem2}
Suppose that fiber transport 
preserves the fiberwise measure $\phi^F \: \dvol_F$ up to
multiplicative constants, 
i.e. for any smooth curve 
$\gamma \: : \: [0,1] \rightarrow B$,
there is a constant $c_\gamma \: > \: 0$ such that
$\rho^* \left( \phi^{F_{\gamma(1)}} \: \dvol_{F_{\gamma(1)}} \right) 
\: = \: c_\gamma \:
\phi^{F_{\gamma(0)}} \: \dvol_{F_{\gamma(0)}}$. \\
1. We have an inequality of functions on $B$ :
\begin{equation}
R^B_\infty \: \ge \: \frac{\int_F (R^M_\infty \: - \: R^F_\infty)
\: \phi^M \: \dvol_F}{\int_F \phi^M \: \dvol_F}.
\end{equation}
2. Suppose in addition that $\phi^M \: = \: 1$. Put $q \: = \: \dim(F)$.
Then
\begin{equation}
R^B_q \: \ge \: \frac{\int_F (R^M \: - \: R^F)
\: \dvol_F}{\vol(F)}.
\end{equation}
\end{theorem}

We now give an application of Theorem \ref{theorem2} to collapsing of
Riemannian manifolds.
A general phenomenon in collapsing is that when the limit
space is smooth, a lower curvature bound on the collapsing manifolds
tends to give a lower curvature bound on the limit space.
This is clearly not true for scalar curvature, as can 
be seen by
taking a product $B \times S^q$, with $q > 1$, and shrinking the
$S^q$-factor. One can get around this particular 
problem by assuming an upper bound
on the scalar curvatures of the collapsing manifolds. The next
result is a scalar curvature analog of \cite[Theorem 3]{Lott (2003)}.

\begin{theorem} \label{theorem3}
Let $\{(M_i, g_i)\}_{i=1}^\infty$ be a sequence of 
$N$-dimensional
connected closed Riemannian manifolds
with sectional curvatures bounded above in absolute value by $\Lambda$
and diameters bounded above by $D$, for some $D, \Lambda \in \R^+$.
Let $(X, \mu)$ be a limit point for 
$\{(M_i, g_i)\}_{i=1}^\infty$ in the measured Gromov-Hausdorff topology.  
Suppose that for some $r \in \R$ and all $i \in \Z^+$, 
$R(M_i, g_i) \: \ge \: r$. 
Suppose that $X$ is an $n$-dimensional closed manifold. 
Put $q \: = \: N \: - \: n$.\\
a. If $q \: = \: 0$
then $X$ has $R \: \ge \: r$.\\
b. If $q \: > \: 0$ then
$X$ has $R_q \: \ge \: r$. Consequently,
$\frac{\int_X R^X \: \dvol_X}{\vol(X)} \: \ge \: r$.
\end{theorem}

In the statement of Theorem \ref{theorem3}, the pair $(g_X, \phi^X)$ is
$W^{2,p}$-regular for all $p \in [1, \infty)$. Hence
$R$ and $R_q$ are well-defined in $\bigcap_{p \in [1, \infty)} L^p(X)$.

The last inequality in Theorem \ref{theorem3} can be compared with the
lower bound in \cite[Corollary 7.49]{Lott-Villani (2004)} on the
mean scalar curvature of a smooth limit space when it arises from a sequence
of collapsing
manifolds with Ricci curvature bounded below by $K \in \R$. 
One could ask whether there is an analog of 
Theorem \ref{theorem3} if one assumes a two-sided bound on the
Ricci curvature of the manifolds $M_i$, as opposed to the sectional curvature.

I thank the UC-Berkeley Mathematics Department 
for its hospitality while this note was written.

\section{Proof of Theorem \ref{theorem2}} \label{Section 3}

We use the notation of \cite[Section 3]{Lott (2003)}, which in turn uses
the notation of \cite[Chapter 9]{Besse (1987)}.
If $X$ is a vector field on $B$, let $\overline{X}$ be its 
horizontal lift to $M$.
Let $N$ be the mean curvature vector field to the
fibers $F$. Let $A$ be the curvature of the horizontal distribution.
Let $T$ be the second fundamental
form tensor of the fibers $F$. Let $\nabla^M$ be the covariant derivative
operator on $M$ and let $\nabla^B$ be the covariant derivative operator on 
$B$.  Define $\check{\delta} N \in C^\infty(M)$ by saying that
if $\{\overline{e}_\alpha\}$ is an orthonormal basis of
$T_{hor}M$ at a point $m$ then
$(\check{\delta} N)(m) \: = \: - \: \sum_\alpha \langle 
\nabla^M_{\overline{e}_\alpha} N, \overline{e}_\alpha \rangle$. 
From \cite[Corollary 9.37]{Besse (1987)}, there is an identity of
functions on $M$ :
\begin{equation} \label{3.1} 
R^M \: = \:
R^B \: + \: R^F \: - \: |A|^2 \: - \: |T|^2 \: - \: |N|^2
\: - \: 2 \: \check{\delta} N.
\end{equation}

Given $f \in C^\infty(M)$,  we write $\nabla^2_M f$ for the Laplacian of
$f$. Define $\nabla^2_{hor} f \in C^\infty(M)$
by saying that $\left( \nabla^2_{hor} f \right)(m)\: = \: \sum_\alpha (\Hess f)
(\overline{e}_\alpha, \overline{e}_\alpha)$.
Define $\nabla^2_F f \in C^\infty(M)$ by saying
that if $m$ is in a fiber $F$ then
$\left( \nabla^2_F f \right)(m)$ is the Laplacian of the function
$f \big|_F$ on $F$, evaluated at $m$. One can check that
\begin{equation} \label{add6}
\nabla^2_M f \: = \: \nabla^2_{hor} f \: + \: \nabla^2_F f \: - \:
\left( \nabla_{hor} f, N \right). 
\end{equation}
  
Given $b \in B$,
let $\{\theta_t\}_{t \in (- \epsilon, 
\epsilon)}$ be the flow of ${X}$ as defined in a neighborhood of
$b$ and for $t$ in some interval $(- \epsilon, 
\epsilon)$. Let $\{\overline{\theta}_t\}_{t \in (- \epsilon, 
\epsilon)}$ be the flow of $\overline{X}$ that
covers $\theta_t$. It sends fibers to fibers diffeomorphically. Hence it 
makes sense to define
${\cal L}_{\overline{X}} \:  \dvol_F$ by 
\begin{equation} \label{3.2}
({\cal L}_{\overline{X}} \:  \dvol_F) \Big|_{F_b} \: = \:
\frac{d}{dt} \Big|_{t = 0} ({\overline{\theta}_t}^* \dvol_F) \Big|_{F_b}.
\end{equation}
With our conventions,
\begin{equation} \label{3.3} 
{\cal L}_{\overline{X}} \:  \dvol_F \: = \: - \: \left( 
\overline{X}, N \right) \:
\dvol_F.
\end{equation}
We have
\begin{equation} \label{3.4} 
\phi^B \: = \: \int_F \phi^M \: \dvol_F. 
\end{equation}
Then
\begin{align} \label{3.5}
X \phi^B \: & = \: {\cal L}_X \phi^B \: = \:
{\cal L}_X \int_F \phi^M \: \dvol_F  = \: 
\int_F {\cal L}_{\overline{X}} \: (\phi^M \dvol_F) \\
& = \: \int_F \left( \overline{X} \phi^M \: - \:
\left( \overline{X}, N \right) \phi^M \right) \: \dvol_F \notag. 
\end{align}
From \cite[(3.7)]{Lott (2003)},
\begin{align} \label{3.7} 
\Hess(\phi^B)(X, X) \: = \:
& \int_F
\left[ \frac{\Hess(\phi^M)(\overline{X},\overline{X})}{\phi^M} \: - \:
\left( \frac{\overline{X} \phi^M}{\phi^M} \right)^2
- \: \left( \overline{X}, \nabla^M_{\overline{X}} N \right)  
\: + \right. \\
& \left. \: \: \: \: \: \:  \: \: \: \: \: \: 
\left( \frac{\overline{X} \phi^M}{\phi^M} \: - \:
\left( \overline{X}, N \right) \right)^2  \right] \phi^M
 \: \dvol_F. \notag
\end{align}
Then
\begin{equation} \label{add1}
\nabla_B^2 \phi^B \: = \:
\int_F
\left[ \frac{\nabla^2_{hor} \phi^M}{\phi^M} \: - \:
\left| \frac{\nabla_{hor} \phi^M}{\phi^M} \right|^2 \:
+ \: \check{\delta} N
\: + \:
\left| \frac{\nabla_{hor} \phi^M}{\phi^M} \: - \:
N \right|^2  \right] \phi^M
 \: \dvol_F.
\end{equation}

Given $b \in B$, put $F_b \: = \: p^{-1}(b)$.
From (\ref{1.1}), (\ref{3.1}) and (\ref{add1}),
\begin{align} \label{add2}
\phi^B(b) \: R^B_\infty(b) \: & = \:
\phi^B(b) \: R^B(b) \: - \: 2 \: \nabla_B^2 \phi^B \: + \:
\frac{|\nabla_B \phi^B|^2}{\phi^B} \\
& = \: 
\int_{F_b} R^B \: \phi^M \: \dvol_F  \: - \: 2 \: \nabla_B^2 \phi^B \: + \:
\frac{|\nabla_B \phi^B|^2}{\phi^B} \notag \\
& = \:
\int_{F_b} 
\left[ \left( R^M \: - \: R^F \: + \: |A|^2 \: + \: |T|^2 \: + \: |N|^2
\right) \: \phi^M \: - \: \right. \notag \\ 
& \left. \: \: \: \: \: \:  \: \: \: \: \: \: 
2 \: \nabla^2_{hor} \phi^M \: + \: 2 \: 
\frac{\left| \nabla_{hor} \phi^M \right|^2}{\phi^M} 
\: - \: 2 \:
\left| \frac{\nabla_{hor} \phi^M}{\phi^M} \: - \:
N \right|^2 \: \phi^M 
 \right] \: \dvol_F 
\: + \notag \\
& \: \: \: \: \: \: \frac{|\nabla_B \phi^B|^2}{\phi^B}. \notag 
\end{align}

We have
\begin{equation} \label{3.10}
{\cal L}_{\overline{X}} (\phi^M \dvol_F) \: = \:
\left( \frac{\overline{X} \phi^M}{\phi^M} \: - \: (\overline{X}, N)
\right) \: \phi^M \: \dvol_F.
\end{equation}
By assumption, $\frac{\overline{X} \phi^M}{\phi^M} \: - \: 
(\overline{X}, N)$ is constant on a fiber $F$.
Let $\{e_\alpha\}$ be an orthonormal basis of
$T_b B$ and let $\{\overline{e}_\alpha\}$
be the horizontal lift of $\{e_\alpha\}$ along $F_b$.  Then
\begin{align}
\frac{\left( e_\alpha \phi^B \right)^2}{\phi^B} \: & = \: 
(\phi^B)^{-1} \: \left( \int_{F_b} 
\left( \frac{\overline{e}_\alpha \phi^M}{\phi^M} \: - \: 
(\overline{e}_\alpha, N) \right) \: \phi^M \dvol_F \right)^2 \\
& = \:
\int_{F_b} 
\left( \frac{\overline{e}_\alpha \phi^M}{\phi^M} \: - \: 
(\overline{e}_\alpha, N) \right)^2 \: \phi^M \dvol_F. \notag
\end{align}
and so
\begin{equation} \label{add3}
\frac{\left| \nabla_B \phi^B \right|^2}{\phi^B} \: = \: 
\int_{F_b} 
\left| \frac{\nabla^M_{hor} \phi^M}{\phi^M} \: - \: 
N \right|^2 \: \phi^M \dvol_F.
\end{equation}
From (\ref{add2}) and (\ref{add3}),
\begin{align} \label{add4}
\phi^B(b) \: R^B_\infty(b) \: & = \:
\int_{F_b} 
\left[ \left( R^M \: - \: R^F \: + \: |A|^2 \: + \: |T|^2 \: + \: |N|^2
\right) \: \phi^M \: - \: \right. \\ 
& \left. \: \: \: \: \: \:  \: \: \: \: \: \: 
2 \: \nabla^2_{hor} \phi^M \: + \: 2 \: 
\frac{\left| \nabla_{hor} \phi^M \right|^2}{\phi^M} 
\: - \: 
\left| \frac{\nabla_{hor} \phi^M}{\phi^M} \: - \:
N \right|^2 \: \phi^M 
 \right] \: \dvol_F \notag \\
& = \:
\int_{F_b} 
\left[ \left( R^M \: - \: R^F \: + \: |A|^2 \: + \: |T|^2
\right) \: \phi^M \: - \: 2 \: \nabla^2_{hor} \phi^M \: + \right. \notag \\ 
& \left. \: \: \: \: \: \:  \: \: \: \: \: \: 
\frac{\left| \nabla_{hor} \phi^M \right|^2}{\phi^M} 
\: + \: 2 \:
\left( \nabla_{hor} \phi^M, N \right)
 \right] \: \dvol_F. \notag 
\end{align}
As 
\begin{equation} \label{add5}
\left| \nabla_M \phi^M \right|^2 \: = \:
\left| \nabla_{hor} \phi^M \right|^2 \: + \: 
\left| \nabla_F \phi^M \right|^2,
\end{equation}
we can combine (\ref{add4}) with (\ref{add6}) and (\ref{add5}) to obtain
\begin{align} \label{add7}
\phi^B(b) \: R^B_\infty(b) \: & = \:
\int_{F_b} 
\left( R^M_\infty \: - \: R^F_\infty \: + \: |A|^2 \: + \: |T|^2
\right) \: \phi^M \: \dvol_F \\
& \ge \: \int_{F_b} 
\left( R^M_\infty \: - \: R^F_\infty 
\right) \: \phi^M \: \dvol_F. \notag
\end{align}
This proves Theorem \ref{theorem2}.1.

Now suppose that $\phi^M \: = \: 1$. From (\ref{add3}), (\ref{add7}) and the
definition of $R_q$ from (\ref{1.2}),
\begin{align}
\phi^B(b) \: R^B_q(b) \: & = \:
\phi^B(b) \: R^B_\infty(b)  \: - \: \frac{1}{q} \:
\frac{|\nabla_B \phi^B|^2}{\phi^B} \\
& = \: 
\int_{F_b} 
\left( R^M \: - \: R^F \: + \: |A|^2 \: + \: |T|^2
 \: - \: \frac{1}{q} \: |N|^2 \right) \: \dvol_F. \notag
\end{align}
As $\left( \overline{e}_\alpha, N \right) \: = \: - \: 
\Tr \left( T \overline{e}_\alpha \right)$, we know that
$\left( T \overline{e}_\alpha, T
\overline{e}_\alpha \right) \:
- \: \frac{1}{q} \: \left( \overline{e}_\alpha, N \right)^2 \: \ge \: 0$.
Summing over $\alpha$ gives
$|T|^2 \: - \: \frac{1}{q} \: |N|^2 \: \ge \: 0$.
This proves Theorem \ref{theorem2}.2.

\section{Proof of Theorem \ref{theorem3}} \label{Section 4}

Let $\{M_i, g_i \}_{i=1}^\infty$ be a sequence as in
the statement of Theorem \ref{theorem3}. We may assume that
$\lim_{i \rightarrow \infty} (M_i, g_i, \dvol_i) \: = \:
(X, \mu)$ in the measured Gromov-Hausdorff topology.
If $q \: = \: 0$ then $X$ is a smooth manifold with a $W^{2,p}$-regular
metric $g_X$ and after taking a subsequence and applying diffeomorphisms, we
may assume that $(M_i, g_i)$ converges to $(X, g_X)$ in the
$W^{2,p}$-topology
(see, for example, 
\cite{Kasue (1989),Peters (1987),Petersen (1997)}). In particular, 
$\lim_{i \rightarrow \infty} R(g_i) \: = \: R(g_X)$ in $L^p(X)$.
Thus $R^X \: \ge \: r$.

Suppose that $q \: > \: 0$.
By saying that $X$ is a manifold, we mean that in the construction of
$X$ as a quotient space $\widehat{X}/O(N)$ \cite{Fukaya (1988)}, the
action of $O(N)$ on the manifold $\widehat{X}$ has a single orbit type.
Then $X$ has the structure of a smooth manifold with a
$W^{2,p}$-regular pair $(g_X, \phi^X)$.

For any $\epsilon \: > \: 0$, we
can apply smoothing results of Abresch and others 
\cite[Theorem 1.12]{Cheeger-Fukaya-Gromov (1992)}
to obtain new metrics $g_i(\epsilon)$ with
\begin{align} \label{4.1}
e^{- \: \epsilon} \: g_i \: \le \: g_i(\epsilon) \: & \le \: 
e^{\epsilon} \: g_i, \\
|\nabla_{g_i} \: - \: \nabla_{g_i(\epsilon)} | \: & \le \: 
\epsilon, \notag \\
|\nabla^k_{g_i(\epsilon)} 
\: \Riem(M_i, g_i(\epsilon))|
\: & \le \: C_{k}(N, \epsilon, \Lambda), \notag
\end{align}
where the constants are uniform.
We can also assume that
$R(M_i, g_i(\epsilon)) \: \ge \: r \: - \: \epsilon$; compare
\cite[Remark 2, p. 51]{Dai-Wei-Ye (1996)} and \cite[Theorem 2.1]{Rong (1996)}.
For small $\epsilon$,
let $B(\epsilon)$ be a Gromov-Hausdorff limit of a subsequence 
of $\{ (M_i, g_i(\epsilon)) \}_{i=1}^\infty$.
We relabel the subsequence as $\{ (M_i, g_i(\epsilon)) \}_{i=1}^\infty$.
From \cite[Proposition 4.9]{Cheeger-Fukaya-Gromov (1992)},
for large $i$, there is a
small $C^2$-perturbation $g_i^\prime(\epsilon)$ of $g_i(\epsilon)$
which is invariant with respect to a $Nil$-structure.
In particular,
we may assume that $R(M_i, g_i^\prime(\epsilon)) \: \ge \: 
r \: - \: 2 \: \epsilon$. Now 
$(M_i, g_i^\prime(\epsilon))$ is the
total space of a Riemannian submersion $M_i \rightarrow B(\epsilon)$ with 
infranil fibers and affine holonomy.
Let $\left( g_i^{B(\epsilon)}, \phi_i^{B(\epsilon)} \right)$ denote the
induced metric and measure on $B(\epsilon)$. As the fiber transport of the
Riemannian submersion preserves the affine-parallel volume forms of the fibers
up to multiplicative constants, and the infranil fibers have nonpositive scalar
curvature, Theorem \ref{theorem2}.2 implies that
$R_q \left(
B(\epsilon), g_i^{B(\epsilon)}, \phi_i^{B(\epsilon)}
\right) \: \ge \: 
r \: - \: 2 \: \epsilon$. 
Varying $i$ and $\epsilon$, we can extract a subsequence of
$\left\{ \left(
B(\epsilon), g_i^{B(\epsilon)}, \phi_i^{B(\epsilon)}
\right) \right\}$ with 
$i \rightarrow
\infty$ and $\epsilon \rightarrow 0$ that converges in
the $W^{2,p}$-topology to $(X, g_X, \phi^X)$. Thus
$R^X_q \: \ge \: r$.

Finally, 
\begin{equation}
R^X_q \: = \: R^X \: - \: 2 \: \nabla^2 \ln(\phi^X) \: - \:
\left( 1 \: + \: \frac{1}{q} \right) \: |\nabla \ln(\phi^X)|^2.
\end{equation}
Then
\begin{equation}
r \: \le \: \frac{\int_X R^X_q \: \dvol_X}{\vol(X)} \: \le \: 
\frac{\int_X R^X \: \dvol_X}{\vol(X)}.
\end{equation}

\end{document}